\documentclass[twocolumn,aps,secnumarabic,nofootinbib,floatfix]{revtex4}
\usepackage{amsmath,amsfonts,amsthm,subfigure,graphics,times,mathptmx}
\newtheorem{theorem}{Theorem}[section]

\newtheorem{lemma}[theorem]{Lemma}
\newtheorem{definition}[theorem]{Definition}
\newtheorem{corollary}[theorem]{Corollary}

\newtheorem{question}[theorem]{Question}
\newtheorem{conjecture}[theorem]{Conjecture}

\makeatletter
\renewcommand{\fnum@figure}[1]{\figurename~\thefigure\ignorespaces}
\renewcommand{\figurename}{Figure}
\makeatother

\newcommand{\fig}[1]{Figure~\ref{#1}}

\newcommand{\R}{\mathbb{R}}
\newcommand{\cross}{\times}

\begin{document}

\title{Quadrisecants of knots and links}
\author{Greg Kuperberg}
\email[Current email: ]{greg@math.ucdavis.edu}
\affiliation{UC Berkeley}

\date{1988}

\begin{abstract}
We show that every non-trivial tame knot or link in $\R^3$ has a
quadrisecant, i.e.  four collinear points.  The quadrisecant must be
topologically non-trivial in a precise sense.  As an application, we
show that a nonsingular, algebraic surface in $\R^3$ which is a knotted
torus must have degree at least eight.
\end{abstract}

\maketitle

\section{Introduction}

An elementary count of degrees of freedom suggests that a randomly-chosen
curve in $\R^3$, if sufficiently complicated, should contain four collinear
points.  One precise interpretation of this intuition is the following
two theorems:

\begin{theorem}[Pannwitz,Morton,Mond] Every non-trivial piecewise linear
or smooth knot in $\R^3$ in general position has four collinear points.  
\end{theorem}

\begin{theorem}[Pannwitz,Morton,Mond] If two smooth or PL circles
$A$ and $B$ in $\R^3$ in general position have a non-zero linking number,
then there is a line in $\R^3$ which intersects $A$, then $B$, then $A$
again, and then $B$ again.
\end{theorem}

These theorems are presented in \cite{MM:quadrisecants} and \cite{Pannwitz:knoten}.  (They are
also mentioned in \cite{BZ:knots}.) Also, the arguments in \cite{Pannwitz:knoten} yield
a lower bound on the number of collinearities and a generalization of
the second theorem to the case of two circles which are linked in the
sense that each represents a non-trivial homotopy class in the
complement of the other.  The main theorem of the present paper is a
different generalization of this result:

\begin{theorem} Every non-trivial tame link in $\R^3$ has four
collinear points. \label{goal}
\end{theorem}

Since the statement of the theorem resembles the statements of theorems
of Pannwitz, Morton, and Mond, we explain the extra cases covered
by our theorem.  By a \emph{non-trivial link} we mean any set of disjoint 
circles embedded in $\R^3$ such that there is no homeomorphism of
$\R^3$ which sends the circles to a flat plane.  The Whitehead
link and the Borromean rings are two examples of non-trivial links
which are not covered by the previous theorems.  A \emph{tame link} is
any set of continuous circles which are collared by solid tori, or
equivalently one which is topologically equivalent to a smooth
link in $\R^3$.  However, a tame link may be very different from
a smooth link geometrically; for example, its Hausdorff
dimension may be greater than 1.  Moreover, the main theorem is not
restricted to links which have any particular transversality properties
or are in general position in any sense.

To eliminate the general position hypothesis, we first prove a stronger
theorem about (smooth) links in general position:  Such a link
has a line which intersects it four times in a topologically
non-trivial way.  Armed with this extra condition, we can use a
limiting argument to pass from links in general position to arbitrary
tame links.

The theorem has an interesting corollary which may be applied to the
topology of real algebraic surfaces.  It is this application which led
the author to the topic of this paper.

\begin{corollary} If an algebraic surface in $\R^3$ contains the
boundary of a knotted solid torus or linked solid tori, the surface
has degree at least 8. \label{degree8}
\end{corollary}

The theorem inspires a definition:

\begin{definition} If $L$ is a link in $\R^3$, a \emph{secant} of $L$ is
a line segment whose endpoints lie in $L$, a \emph{trisecant} of $L$ is
a secant of $L$ and a point $p$, the \emph{middle point}, which lies in
both $L$ and the interior of the secant, and a \emph{quadrisecant} is a
secant with two middle points.  \end{definition}

To be precise, a quadrisecant is a pair of distinct trisecants with the
same underlying line segment.  A \emph{degenerate secant} is a single
point.  The set of secants has a natural topology, as does the set of
trisecants: for a sequence of trisecants to converge we insist that the
middle points converge as well.

As motivation for the main theorem, we present a simple proof of a
weaker result:

\begin{theorem} Every non-trivial smooth knot $K$ in $\R^3$ has a trisecant.
\end{theorem}

\begin{proof} Suppose that there exists a point $p$ in $K$ such that no
points $q$ and $r$ in $K$ are collinear with $p$.  Then the union of the
chords $\overline{pq}$ for all $q$ in $K$ is evidently a smooth
embedded disk with boundary $K$, which renders $K$ trivial.
\end{proof}

This proof illustrates the central idea in the proof of
the main theorem.

I would like to thank my advisor, Andrew Casson, for encouragement and
helpful comments.

\section{General position}

There is a general theory of general position, presented in a paper by Wall
\cite{Wall:generic} and used in \cite{MM:quadrisecants}.  We review the
elements of this theory needed here:

\begin{definition} If $X$ is a topological space with a measure, a
property $P$ of members of $X$ is \emph{generic} if it is true on a set
with full measure, and a member of $X$ is in \emph{general position}
with respect to $P$ if it satisfies $P$.  A member of $X$ is in
\emph{general position} if it is in general position with respect
to all applicable generic properties mentioned in this paper.
\end{definition}

Usually $X$ is a space of functions.  We define a \emph{polynomial
function} from the unit circle $S^1$ in $\R^2$ to $\R^3$ to be a
function which is given by polynomials of some degree $d$ in the
standard coordinates in $\R^2$.  The set of all such functions forms a
finite-dimensional vector space $P_d$, and we consider all generic
properties relative to $P_d$ for some $d>0$ with the usual Cartesian
topology and measure.  A function $K:S^1 \to \R^3$ is a \emph{knot} if
it is injective, and since this is a generic property, we are justified
in referring to polynomial functions in general position as polynomial
knots.

More generally, we may define $k(S^1)$ to be the disjoint union of $k$
unit circles, consider the vector space $P_{d,k}$ of $k$-tuples
of polynomial functions, and define a \emph{link} to be an injective
function from $k(S^1) \to \R^3$ for some $k$.

The concept of a polynomial link is not an essential ingredient in this
paper, but the following lemmas, whose proofs are easy, make it a useful one:

\begin{lemma} Given an arbitrary smooth function $f:k(S^1) \to \R^3$,
there is a sequence of polynomial links (of varying degree) whose
values and first derivatives converge uniformly to those of $f$.  We
may choose the sequence to be in general position.  \end{lemma}

\begin{lemma} A property $P$ of members of a finite-dimensional vector space
is a \emph{polynomial property} or an \emph{algebraically generic} property
if there exists some non-trivial polynomial $p$ on the 
vector space such that $P$ is true at all points for which $p$ is
non-zero.  All polynomial properties are generic.  
\end{lemma}

If $L$ is a polynomial link with $k$ components, we define a projection
function $\pi_L:k(S^1) \cross k(S^1) - \Delta \to S^2$, where $\Delta$
is the diagonal, by: $$\pi_L(a,b) = (L(a) - L(b))/|L(a) - L(b)|.$$ We
view $\pi_L$ as a family of maps $\pi_L(\cdot,b)$ parameterized by the
second variable.

The main result of this section is the following lemma.  Neither the lemma nor
the proof have more mathematical content than equivalent lemmas in
\cite{MM:quadrisecants} and \cite{Pannwitz:knoten}, and the key idea is
originally due to Reidemeister \cite{Reidemeister:knotentheorie}, so the proof
here is sketched to some extent.

\begin{figure}
\begin{center}
\scalebox{.25}{\includegraphics{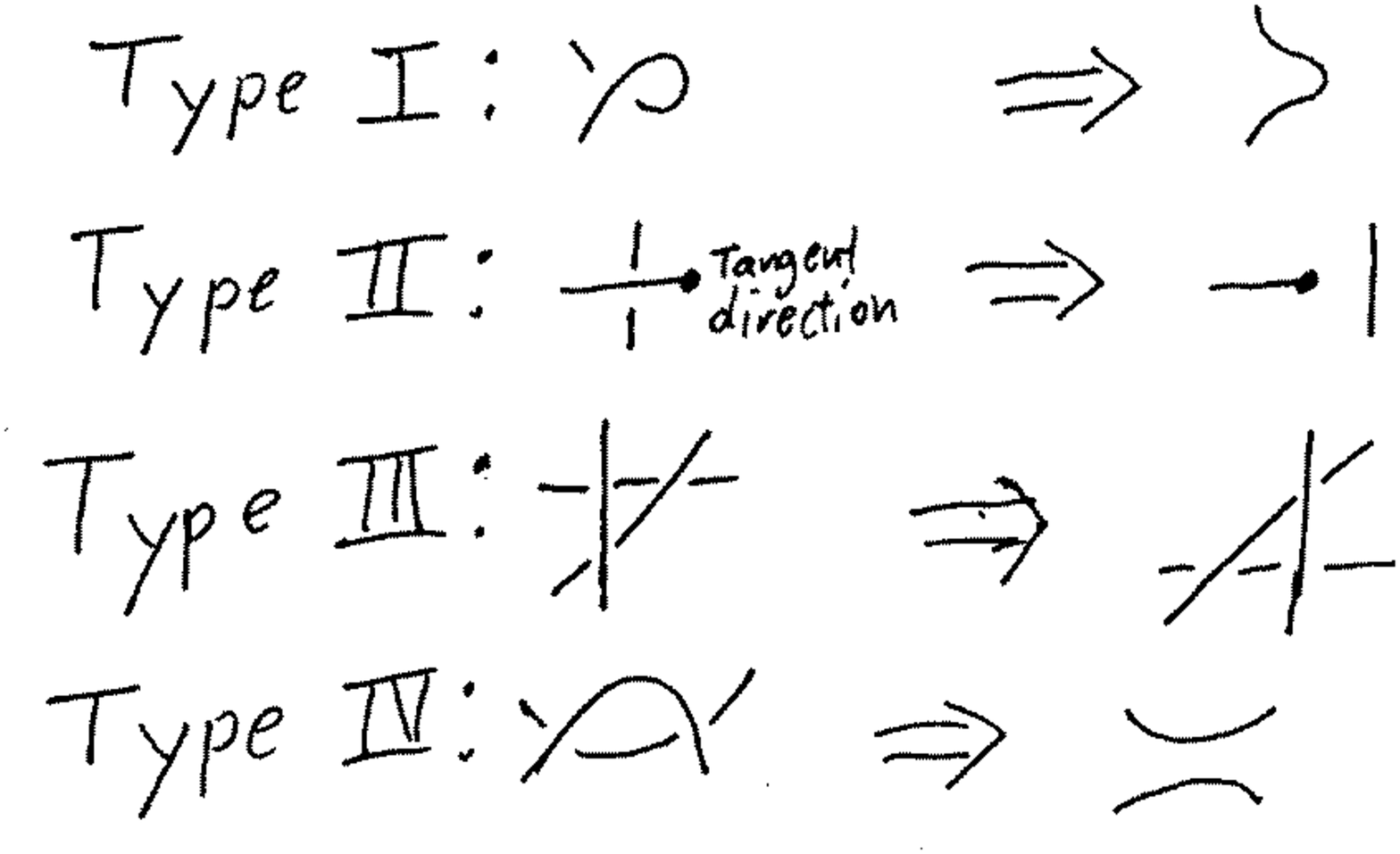}}
\caption{\label{f:one}}
\end{center}
\end{figure}

\begin{lemma} With $L$ and $\pi_L$ defined as above, it is a polynomial
property for $L$ to be a \emph{smooth embedding}, i.e. its derivative
does not vanish anywhere.  It is also a polynomial property of $L$ for
there to exist a finite set of points of $k(S^1)$, called the set of
\emph{special points}, whose complement is the set of \emph{generic
points}, such that for a generic point $a$ and a special point $b$:
\label{moves}
\begin{description}
\item[I.] $\pi_L(\cdot,a)$ is a smooth immersion of a 1-manifold with ends,
where the ends correspond to the tangent directions of $L$ at $a$.
\item[II.] $\pi_L(\cdot,a)$ does not pass through the two tangent directions.
\item[III.] $\pi_L(\cdot,a)$ is everywhere one-to-one or two-to-one.
\item[IV.] If $\pi_L(\cdot,a)$ is two-to-one at a point of $S^2$, it is
self-transverse at that point.
\item[V.] $\pi_L(\cdot,b)$ has all of the previous properties at all but one
point of $S^2$ and has three of the previous properties at the remaining
point $p$.  In this case, as $a$ varies from one side of $b$ to the other,
the structure of $\pi_L(\cdot,a)$ near $p$ is characterized by one of the
corresponding diagrams in \fig{f:one}.
\end{description}
\end{lemma}

\begin{proof} We define the \emph{algebraic dimension} of a subset $S$
of a vector space $V$ to be the Krull dimension of the ring of
polynomial functions restricted to $S$.  (The Krull dimension of a
commutative ring is the maximum length of an ascending chain of prime
ideals \cite{Hartshorne:gtm}.)  We will need two basic facts about algebraic
dimension:  The algebraic dimension image of a set $S$ under a
projection (or more generally a polynomial map) is less than or equal
to the algebraic dimension of $S$, and the complement of a set of
algebraic codimension 1 or more is a polynomial property.  In
the following discussion we will also use codimension to mean
the difference of the dimension of a pair of sets.

For simplicity, we consider only the case of knots.  Observe
that in the vector space of ordered quadruples of points in $\R^3$,
the subset for which the four points are collinear has algebraic
codimension 4.  Given four points $a$, $b$, $c$, and $d$ on
the unit circle, the space of knots $K$ of degree $d$ (for $d \ge 2$)
projects onto the space of quadruples of points in $\R^3$.
Therefore the set of knots $K$ of degree $d$ for which $K(a)$, $K(b)$,
$K(c)$, and $K(d)$ are collinear has codimension 4 as well, 
as does the analogous set in the space of quintuples $(K,a,b,c,d)$,
where $a$, $b$, $c$, and $d$ are four distinct points on the circle.
By projection, the set of pairs $(K,a)$ for which there exists
$b$, $c$, and $d$ such that $K(a)$, $K(b)$, $K(c)$, and $K(d)$
are collinear has codimension at least 1.  Except for an algebraic
subset of the set of knots, the set of $a$ for a knot $K$ for
which $b$, $c$, and $d$ can be found with this property is polynomial,
i.e. finite.  Such a $b$, $c$, and $d$ would have to exist
in order for $\pi_L(\cdot,a)$ to be three-to-one.  Thus, part III
of the lemma is proved for knots.

The rest of the lemma can proved in the same fashion, namely by keeping
track of the codimension of certain sets.  Informally, a set of
algebraic codimension $n$ can be called an $n$-fold coincidence.  Parts
I and II of the lemma hold because, given points $a$ and $b$ on a link
$L$, it would take a 2-fold coincidence for the tangent to $L$ at $b$
to contain $a$, and allowing $a$ to vary, it would take a 1-fold
coincidence in the choice of $a$, or allowing $a$ to vary, a 1-fold
coincidence for the choice of $b$.  Part IV of the lemma holds because,
given $a$, $b$, and $c$ on a link $L$, it would take a 3-fold
coincidence for $a$, $b$, and $c$ to be collinear and for the tangent
lines at $b$ and $c$ to be coplanar.

For part V, the case when condition IV of the lemma fails typifies the
method of proof.  Informally, at a special point $a$ for which
$\pi_L(\cdot,a)$ is somewhere three-to-1, three arms of the projection
of the link meet at a point and it would take a coincidence for there
to be a fourth arm at the point or for two of the arms to have the same
slope.  Near $a$ the front two arms cross at a point and it would take
a coincidence for that crossing to travel parallel to the third arm
instead of passing through it.

Geometrically, it would take a 6-fold coincidence for five given points
on a link $L$ to be collinear, and it would take a 5-fold coincidence
for four given points on $L$ to be collinear and for two of the tangent
lines to be coplanar.  So in either case it would take a 1-fold
coincidence in the choice of $L$ for such a set of points to exist.
Finally, consider collinear four points $a$, $b$, $c$, and $d$ on $L$
and let $l_a$, $l_b$, $l_c$, and $l_d$ be the tangent lines at these
points.  The set of lines that intersect $l_a$, $l_b$, and $l_c$ sweeps
out a surface, and it is a 1-fold coincidence in the choice of $l_d$
for it to be tangent to that surface.  If it is not tangent, then
$\pi_L(\cdot,a)$ will look as it does in case IV of \fig{f:one}.
\end{proof}

\section{Knots in general position}

The arguments in this section follow that of \cite{Pannwitz:knoten} and
\cite{MM:quadrisecants}. The only new feature is the notion of topological
non-trivial quadrisecants, which we will need to generalize the main theorem to
arbitrary knots.

We begin with a simple lemma and a definition:

\begin{lemma} Let $C$ be a compact set in $\R^n$.  Then not every point
of $C$ lies between two other points of $C$. \label{compact}
\end{lemma}
\begin{proof} If $p$ is any point in $\R^n$, then a point $q \in C$
which is farthest from $p$ has this property, because if $q$ lay between
two other points, one of them would be still farther away.
\end{proof}

\begin{definition} A secant of a link $L$ with no extra interior intersections
with $L$ is \emph{topologically trivial} if its endpoints lie on the same
component of $L$, and if it, together with one of the two arcs of this
component, bounds a disk whose interior does not intersect $L$.  The disk may
intersect itself and the secant.  A quadrisecant $\overline{ad}$ with middle
points $b$ and $c$ is topologically trivial if any of the secants
$\overline{ab}$, $\overline{bc}$, and $\overline{cd}$ are.  Similarly for a
trisecant.
\end{definition}

\begin{lemma} A knot in general position has a topologically non-trivial
quadrisecant.
\end{lemma}

\begin{figure}[ht]
\begin{center}
\subfigure[]{\scalebox{.35}{\includegraphics{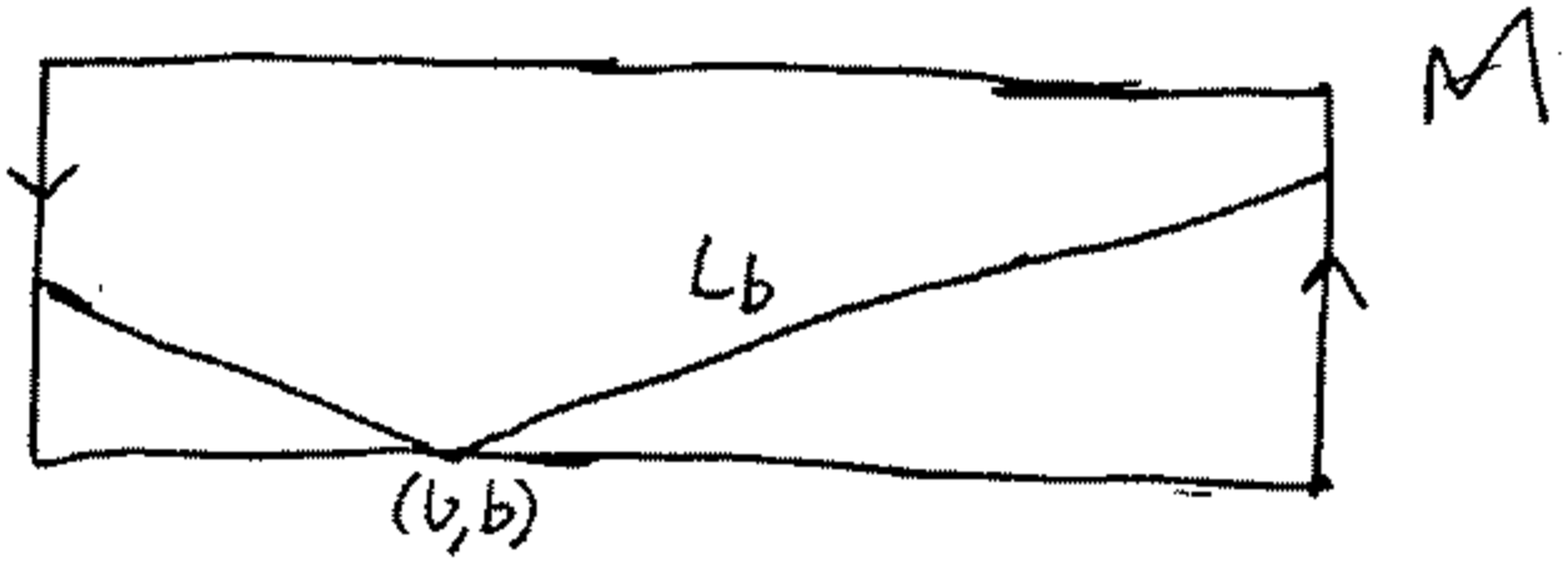}}}
\subfigure[]{\scalebox{.35}{\includegraphics{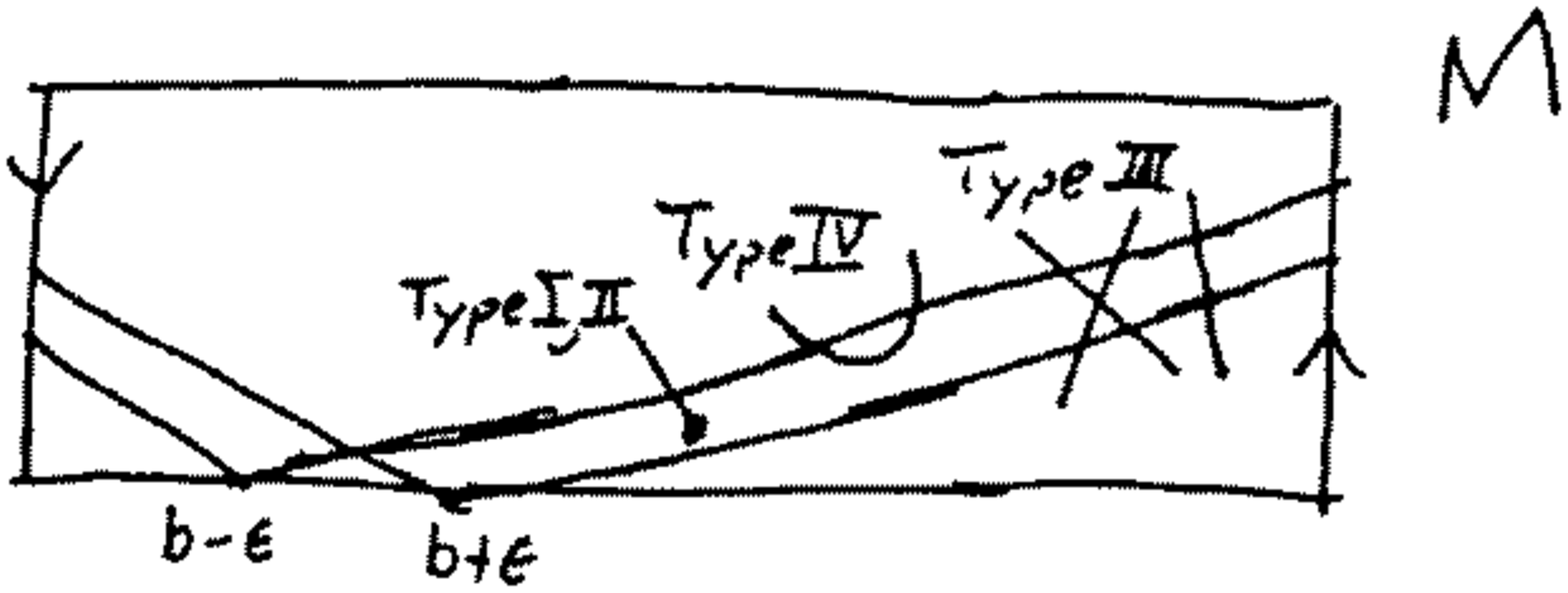}}}
\caption{\label{f:two}}
\end{center}
\end{figure}

\begin{proof}
Let $K$ be a polynomial knot in general position.  Let $M$ be the set of
unordered pairs of points of $S^1$, or equivalently the set of secants of $K$. 
$M$ is topologically a M\"{o}bius strip.  We define $O$ to be the subset $M$
consisting of those pairs of points $(a,b)$ with the property that at least one
point of $K$ lies between $K(a)$ and $K(b)$.  Lemma~\ref{moves} has
implications about the local structure of $O$.  For fixed $a$, the set $L_a$ of
all $(a,b)$ in $M$ is a line segment which wraps around $M$ as in \fig{f:two}(a). 
The intersection $O \cap L_b$ is a finite set.  If $b$ is a generic point, the
topology of $\pi_K(\cdot,a)$, and therefore the topology of $O \cap L_a$,
cannot change as we vary $a$ slightly.  But if $a$ is a special point, the
topology of $O \cap L_a$ changes as illustrated in \fig{f:two}(b).  For example, if
$a$ is a special point at which condition IV of Lemma~\ref{moves} for
$\pi_K(\cdot,a)$ fails, then there exist three points $b$, $c$, and $d$ so that
$a$, $b$, $c$, and $d$, in that order, make a quadrisecant of $K$.  The
trisecants $a$, $c$, $d$ and $a$, $b$, $d$ represent the same point of $O$, and
if condition V of Lemma~\ref{moves} holds, they represent arms of $O$ that
cross. Meanwhile the trisecant $a$, $b$, $c$ represents a point of $O$ that
lies elsewhere along $L_a$.

It follows that $O$ is the image of a self-transverse smooth immersion
of a 1-manifold, and a self-crossing corresponds to a quadrisecant.
If $C$ is a curve of points in $O$ which does not ``make turns''
at the self-crossings, then $C$ is a continuous curve of trisecants.

The significance of $O$ is that it is an obstruction to the following
construction:  Recall that on a M\"{o}bius strip, there are two kinds
of properly embedded arcs, non-separating arcs and separating arcs.
Suppose that $A$ is a non-separating arc of $M$ which avoids $O$.  Then
$A$ corresponds to a family of secants whose interiors do not intersect
$K$.  This family of secants induces a map $D$ from the unit disk to
$\R^3$ whose boundary is $K$ and whose interior does not intersect $K$.
By Dehn's lemma, $K$ is trivial.

Suppose that $K$ has no quadrisecants.  Then $O$ is an embedded 
1-manifold.  By elementary homology theory, if $O$ obstructs all
non-separating arcs, there is a circular component $C$ of $O$ which 
winds around $M$ either one or two times.  The curve $C$ is a continuous family
of trisecants.  We consider the corresponding families of points $\{a,b\}_t$
and $m_t$, with $t \in S^1$, such that $K(m_t)$ lies between
$K(a_t)$ and $K(b_t)$.  If $C$ winds once around $M$, the endpoints travel
half way around $S^1$ and then switch places, and since $m_t$ is trapped
between them, it must jump discontinuously, a contradiction.  If $C$
winds twice around $M$, the endpoints each wind once around $S^1$, and
therefore so does $m_t$.  Thus, every point of $K$ lies between two other
points, which contradicts Lemma~\ref{compact}.

Topological non-triviality is achieved by a modification of this construction.
Let $O'$ be the subset of $O$ consisting of topologically non-trivial
trisecants and quadrisecants which are non-trivial at the middle secant.
Observe that $O'$ is also the image of a smooth immersion:  If $O'$ contains
a self-intersection point of $O$ but does not contain all four arms
of the self-intersection, then it must contain exactly two arms, and
they must be opposite rather than adjacent.   In this case the
self-intersection point is a quadrisecant which is topologically trivial
on one side.  Therefore if $O'$ has a self-crossing, it corresponds to
a topologically non-trivial quadrisecant.  If there are no such quadrisecants,
$O'$ must also have a circular component $C$ with all of the properties
mentioned above, provided that $O'$ is also an obstruction to all
non-separating arcs $A$.

\begin{figure*}[ht]
\begin{center}
\subfigure[]{\scalebox{.35}{\includegraphics{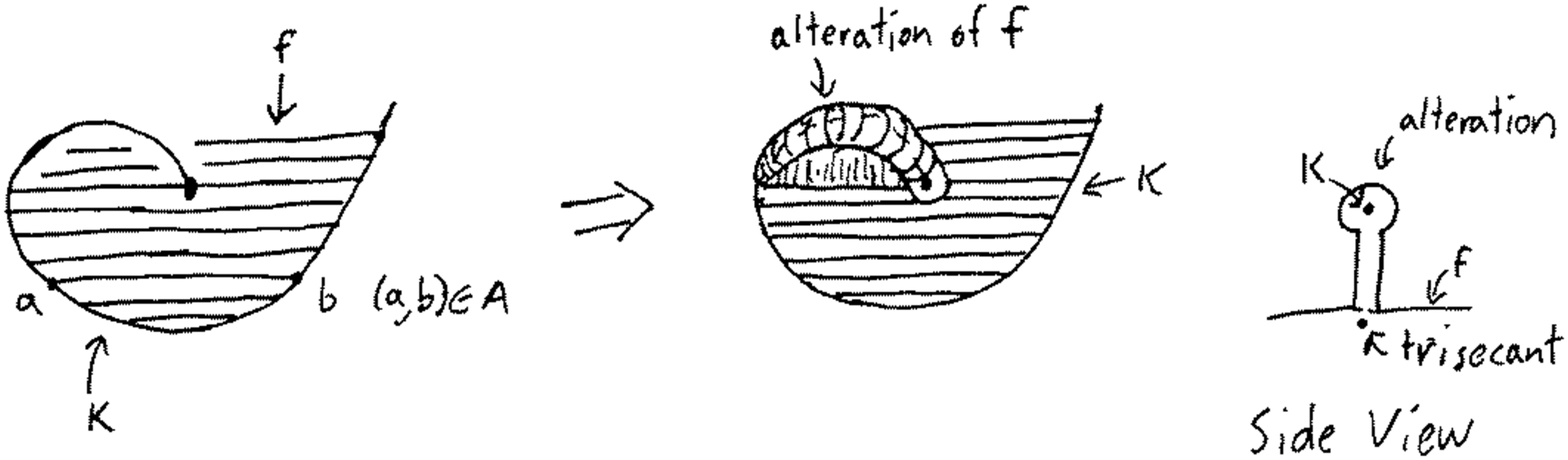}}}
\subfigure[]{\scalebox{.35}{\includegraphics{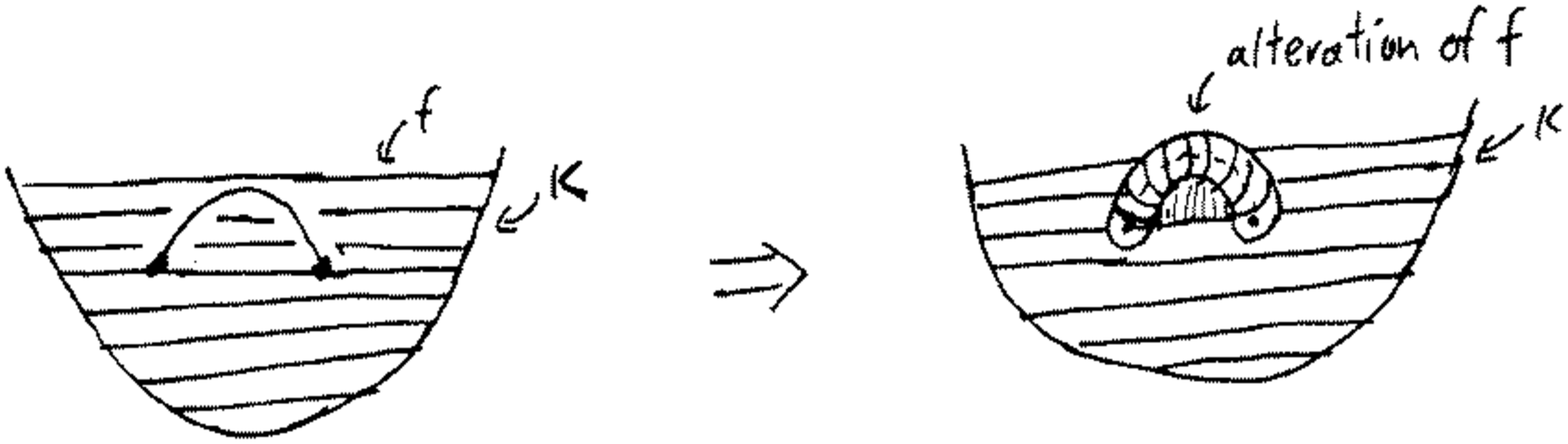}}}
\caption{\label{f:three}}
\end{center}
\end{figure*}

Let $A$ be a non-separating arc which avoids $O'$.  We may choose $A$
to be transverse to $O$.  As before, we construct the disk $D_A$ from
the secants of $A$, but this time $D_A$ does not avoid $K$.  Consider
a point where $A$ intersects a topologically trivial trisecant $T$.
By hypothesis there exists a disk $D_T$ which bounds a secant of $T$
and an arc of $K$.  Using $D_T$ and a tubular neighborhood of $K$, we
may alter $D_A$ to obtain a disk $D'_A$ which avoids $K$ in the vicinity
of $T$, according to \fig{f:three}(a).  We may similarly modify $D_A$ in the
vicinity of a quadrisecant $Q$ which is topologically trivial in the
middle, as in \fig{f:three}(b).  In this fashion we obtain a disk $D$ whose
interior avoids $K$ as before, and Dehn's lemma applies.
\end{proof}

\section{Links in general position}

The result of this section is a completion of analogous results in
\cite{MM:quadrisecants} and \cite{Pannwitz:knoten}.  The arguments there roughly
correspond to the $omega_1 \ne 0$ case of the proof, although the
argument in \cite{Pannwitz:knoten} is somewhat more general than this special
case.

\begin{lemma} Every non-trivial link $L$ in general position
has a topologically non-trivial quadrisecant.
\end{lemma}
\begin{proof} We may assume without loss of generality that no
component of $L$ bounds a disk whose interior avoids $L$.

Let $K$ be a component of $L$.  Let $M_K$ be the M\"{o}bius strip of
secants of $K$, and let $O'_K$ be the corresponding set of
topologically non-trivial trisecants and quadrisecants which are
non-trivial in the middle.  As before, $O'_K$ must be an obstruction to
non-separating arcs $A$, and we obtain a circle $C$ which winds around
$M$.  If the middle points of $C$ also lie on $K$, we may apply the
proof of the previous lemma.  But the middle points may lie on some
other component $H$ of $L$.  In this case, the secants of $C$ induce a
map $f$ from a surface $E$ to $\R^3$, where $E$ is either an annulus or
a M\"{o}bius strip, depending on whether $C$ winds once or twice around
$K$.  We may choose $f$ so that the median of $E$ maps to the middle
points of the trisecants of $C$.

The set of lines $l$ perpendicular to $H$ at a given point $p$ is
homeomorphic to a circle, and the corresponding set $T$ of all ordered
pairs $(l,p)$ is homeomorphic to a torus.  We may orthogonally project
each trisecant $t \in C$ to a line perpendicular to $H$, i.e. a
member of $T$, thereby obtaining a map $f$ from $C$ to $T$.  Since $C$ is a
circle, this map has an ordered pair of winding numbers
$(\omega_1,\omega_2)$ which are well-defined up to an orientation of $C$.
There are three cases to consider, depending on the values of the winding
numbers.

Suppose that $\omega_1 = \omega_2 = 0$.  We construct a disk whose
boundary is $K$ and whose interior avoids $L$.  The map $f$ intersects $H$ at
the median, but it may also intersect $K$ at some other points, because
$C$ may include some quadrisecants which are topologically trivial on one
side.  In this case we can modify $f$ according the prescription in
\fig{f:three}(a) to obtain a map $f'$ which avoids $K$ in the interior and
which agrees with $f$ in a neighborhood of the median.  Since both
winding numbers are zero, we may now homotop $f'$ in a neighborhood
of $H$ to obtain a map $f''$ which avoids $H$ and is constant on
the median of $E$.  Finally, we identify the median of $E$ to a point
to obtain a space $E'$ and a map $f'''$.  If $E$ is a M\"{o}bius
strip, $E'$ is a disk, but if $E$ is an annulus, $E'$ is two disks
identified at a point.  Either way, we obtain the desired spanning
disk, which we may convert to an embedded disk by Dehn's Lemma.

Suppose instead that $\omega_2 = 0$ but $\omega_1 \ne 0$.  Then we
extend $E$ to a line bundle $E'$ and extend $f$ linearly to a map
$f':E' \to \R^3$.  We can homotop $f'$ in a neighborhood of $H$ without
changing its values in $E'\backslash E$ to a map $f''$ which has
constant value $p$ on the zero section of $E'$, but we cannot make
$f''$ avoid $H$.  Let $p \in H$.  As before, we identify the zero
section of $E'$ to a point and obtain a space $E''$, and
correspondingly alter $f''$ to obtain a map $f''':E'' \to \R^3$.  This
time the intersection number between $H$ and $f'''$ at $p$ is
$\omega_1$.  But since $f'''$ is a closed map from the pseudo-manifold
$E''$ to $\R^3$, it induces a well-defined homology class in the
infinite homology of $\R^3$.  $H$ induces another such homology class,
and by elementary homology theory, the total intersection number
between $f'''$ and $H$ must be zero.  The map $f'''$ must intersect $H$
at another point, and therefore $f''$ does also.  Suppose that $f''(x)$
is this point, with $x \in E'$. The point $x$ cannot be in $E$,
therefore $f'(x) = f''(x)$.  Since $f'$ is linear on the fibers of
$E'$, the image under $f'$ of the fiber containing $x$ yields
a quadrisecant $Q$.  The quadrisecant $Q$ is necessarily topologically
non-trivial, because if the intersection points of $Q$ are labeled in
order as $a,b,c,$ and $x$, then $b,x \in H$ and $a,c \in K$.

The only remaining possibility is that $\omega_2 \ne 0$.  In this
case, every point of $H$ lies between two points of $K$.  We may
repeat the whole argument with each component of $L$ playing the
role of $K$, thereby obtaining a function $f$ from components
of $L$ to components of $L$ such that every point of $f(K)$
lies between two points of $K$.  The map $f$ must have at least one
circular orbit, and we may set $C$ to be the set in $\R^3$
which is the union of all components of $L$ in this orbit.  
Evidently, $C$ is a compact set and every point of $C$ lies
between two other points of $C$, a contradiction by
Lemma~\ref{compact}.
\end{proof}

\section{Arbitrary tame knots and links}

\begin{definition}  A link $L$ in $\R^3$ is \emph{tame} if there 
exists a homeomorphism $h$ of $\R^3$ which carries $L$ to a polynomial
link, or equivalently a piecewise linear or smooth link.
\end{definition}

\begin{lemma} If $L$ is a tame link, there exists a homeomorphism
$h$ of $\R^3$ which maps $L$ to a smooth link with $h$ smooth on
$\R^3 - L$. \label{smoothing}
\end{lemma}
\begin{proof} Let $K$ be a tame knot and let $h$ be an arbitrary homeomorphism
such that $h(K)$ is smooth.  Using a tubular neighborhood of $h(K)$, we can
choose $T_i$, with $i\ge 1$, to be a sequence of nested, parallel tori
converging to $h(K)$.  Let $T'_i = h^{-1}(T_i)$. By the theory of
triangulations and smoothings of 3-manifolds (see \cite[p.~217]{Moise:gtm}),
there exists a sequence of smooth tori $T''_i$, with each $T''_i$ lying between
$T'_i$ and $T'_{i+1}$, and a sequence of diffeomorphisms $h'_i:T''_i \to T_i$
such that ${h'}^{-1}_i$ and $h^{-1}|_{T_i}$ are isotopic as maps from $T_i$ to
$\R^3 - K$. Furthermore, we can arrange that the distance between ${h'}^{-1}_i$
and $h^{-1}$ goes to zero as $i \to \infty$.  By the isotopy condition, the
$h'_i$'s may be extended smoothly to the each region between $T'_i$ and
$T'_{i+1}$ and the region outside $T_1$ to obtain a diffeomorphism $h':\R^3-K
\to \R^3 - h(K)$.  Because of the distance condition, we can continuously
extend $h'$ to $K$ by setting it equal to $h$ on $K$. This continuous extension
is the desired map.

The proof in the case of links is similar.
\end{proof}

We are now in a position to prove Theorem~\ref{goal}.  In fact,
we can prove something slightly stronger:

\begin{theorem} If $L$ is a non-trivial tame link in $\R^3$, then $L$
has a quadrisecant, none of whose component secants are subsets of $L$.
\end{theorem}

\begin{proof} Let $L$ be a non-trivial, tame link and let $h$ be a
homeomorphism given by Lemma~\ref{smoothing}.  Let $N$ be a tubular
neighborhood of $h(L)$, let $N'$ be the normal bundle of $h(L)$, and
choose a diffeomorphism $n:N \to N'$.  Consider a sequence of links
$L_i$ such that $h(L_i)$ is disjoint from $L$ and $n(h(L_i))$ is a
smooth section.  Choose the sequence so that $h(L_i)$ converges
smoothly to $h(L)$, i.e. $n(h(L_i))$ converges smoothly to the zero
section.  Since $h$ is a diffeomorphism outside of $L$, we may choose
each $L_i$ to be a polynomial link in general position.

By hypothesis, each $L_i$ has the same isotopy type as $L$, and in
particular each $L_i$ is non-trivial.  Therefore each $L_i$ has a
topologically non-trivial quadrisecant $Q_i$.  By compactness,
$\{Q_i\}$ has a convergent subsequence in the space of line segments in
$\R^3$; we may suppose without loss of generality that the original
sequence converges.  The resulting limit is a secant of $L$.  We must
show that the endpoints and middle points of the quadrisecants do not
converge together.

For each $i$, let $S_i$ be a topologically non-trivial secant of $L_i$
and suppose that the $S_i$'s converge to a point $p$ on $L$.  Let $B$
be a round, open ball in $N'$ centered at $n(h(p))$.  Then there exists
an $i$ such that $S_i$ and an arc $A$ of $L_i$ with the same endpoints
as $S_i$ are both contained in $h^{-1}(n^{-1}(B))$.  For each point $s
\in n(h(S_i))$, we consider the line segment from $s$ to $t$, where $t$
is the point in $n(h(L_i))$ which lies in the same fiber of $N'$ as
$s$, as illustrated in \fig{f:four}.  Since $n(h(L_i))$ is a section, $t$
is unique.  The union of these line segments is the image of a spanning
disk of $n(h(A \cup S_i))$ which does not intersect $n(h(S_i))$.
Therefore $S_i$ is topologically trivial, a contradiction.  

\begin{figure}
\begin{center}
\scalebox{.3}{\includegraphics{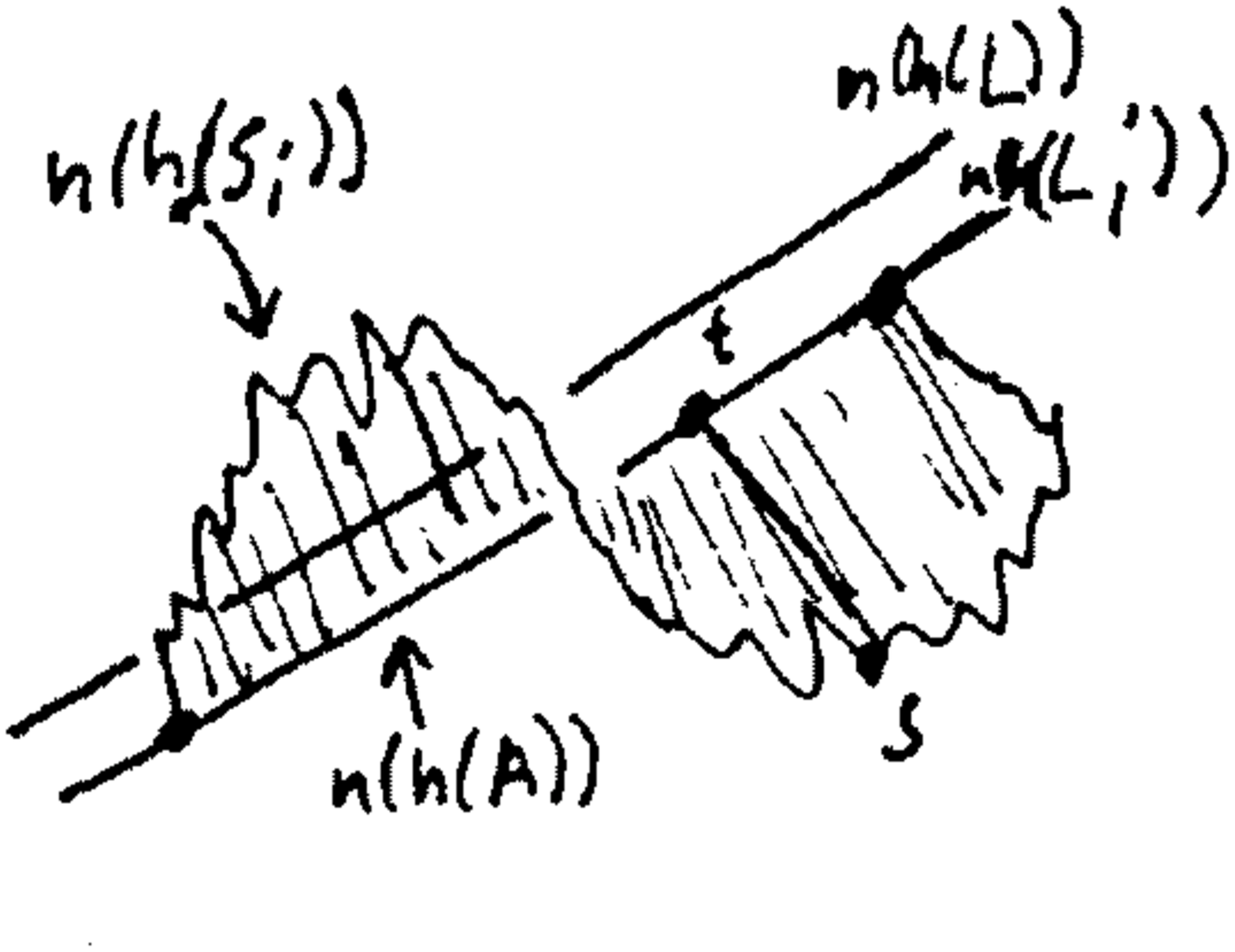}}
\caption{\label{f:four}}
\end{center}
\end{figure}

The proof that the limit of the $S_i$'s is not a subset of $L$ is similar.
\end{proof}

Corollary~\ref{degree8} follows from this theorem:

\begin{proof} Let $\{T_i\}$ be a non-trivially linked collection of 
solid tori.  For each $i$ and each $n>0$, let $l_{i,n}$ be the shortest
non-contractible loop in $T_i$ which is homotopically $n$ times the
core of $T_i$.  Let $l_i$ be a shortest member of the set $\{l_{i,n}\}$.
If we let $D$ and $D'$ be two disjoint, non-separating disks in $T_i$
for some $i$, then we see that the length of $l_{i,n}$ is bounded below
by $n$ times the distance between $D$ and $D'$.  Therefore $l_i$ exists,
although it may not be unique.

Suppose that for some $a,b \in S^1$, $l_i(a) = l_i(b)$.  Then we can divide
$l_i$ into two loops from $l_i(a)$ to itself.  At least one of these loops must
be non-contractible and both loops are shorter, which is a contradiction. 
Thus, each $l_i$ is an embedding.  If we let $L$ be the union of the images of
the $l_i$'s, then $L$ is a satellite link of the $T_i$'s.  By a theorem in knot
theory \cite[p.~113]{Rolfsen:knots}, $L$ must be a non-trivial link if the
$T_i$'s are.

Since a geodesic in a smooth manifold with smooth boundary must be
$C^1$ (see \cite{AA:geodesics}; a proof was also suggested to the author by
Tom Ilmanen), $L$ must be a tame link.  By the preceding theorem, $L$
must have a quadrisecant $Q$ such that no component secant of $Q$ is
contained in $L$.  Suppose that a component secant $S$ of $Q$ were
contained entirely inside some $T_i$.  Let $p$ be a path which goes
from one endpoint of $S$ to the other.  Then we can divide $l_i$ into
two paths $q_1$ and $q_2$ to make two loops $q_1p$ and $q_2p$ whose
composition is homotopic to $l_i$.  At least one of these loops must be
non-contractible, therefore they cannot both be shorter.  Therefore
each component secant of $Q$ must have one point which lies outside the
$T_i$'s.

Finally, suppose that $P(x,y,z)$ is a non-trivial polynomial whose zero set
contains $\partial T_i$ for all $i$.  Then the restriction of $P$ to the line
containing $Q$ must be non-trivial and must have at least 8 real roots,
counting multiplicity.  Therefore $P$ has degree at least 8.
\end{proof}

\begin{figure*}[ht]
\begin{center}
\scalebox{.35}{\includegraphics{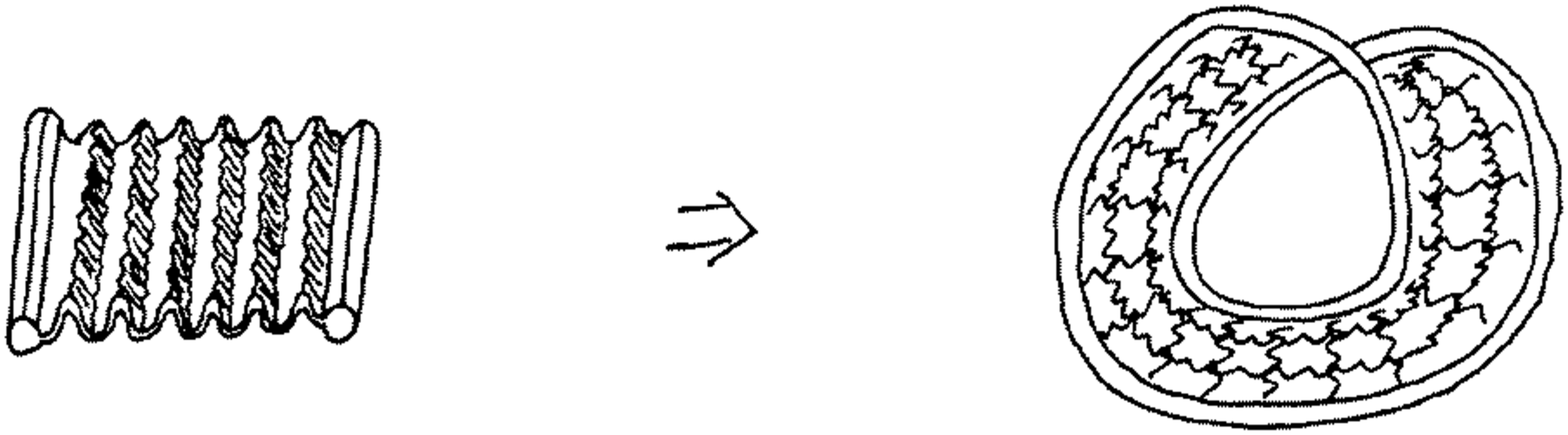}}
\caption{\label{f:five}}
\end{center}
\end{figure*}

The author once believed that the loop in a solid torus which is the
shortest non-zero multiple of the core is necessarily homotopic to
the core.  However, this is false by an example of Doug Jungreis.
We can consider the region $S$ in $\R^3$ which consists of the set of 
points $(x,y,z)$ such that:
$$|x - \sin(L_1^2y)/L_1 - \sin(L_2^2z)/L_2| < \epsilon,$$
where $L_1$ is very large, $L_2$ is much larger still, and $\epsilon$
is much smaller than $1/L_2$.  The region $S$ could be described as a
corrugated sheet, and it has the property that if $a,b \in S$ and the
straight-line distance from $a$ to $b$ is greater than 1, then this
distance is much less than the length of the shortest path in $S$ from
$a$ to $b$.  If $M$ is a smooth M\"obius strip in $\R^3$ whose
tangent plane varies slowly, we can approximate $M$ with a solid
torus $T$ which is topologically a tubular neighborhood of $M$
but which is geometrically quite different:  $T$ is the union of
a thick tube centered around the boundary of $M$ and a corrugated sheet
which approximates the interior of $M$, as shown in \fig{f:five}.
Clearly the shortest non-trivial loop in $T$ stays close to the
boundary of $M$ and is therefore homotopically twice the core.

It is easy to show that the bound in Corollary~\ref{degree8} is the 
best possible:  If we choose two numbers $r_1 > r_2$, then the surface
given by:
\[ (x^2+y^2+z^2-r_1^2-r_2^2)^2 - 4(x^2+y^2)r_1^2 = 0 \]
is a torus.  If $r_1 > 2r_2$, we can multiply two such surfaces together
to obtain two linked tori.

\section{Questions open to the author}

The most serious shortcoming of Corollary~\ref{degree8} is the fact
that it only applies to closed surfaces in $\R^3$, while the usual
context for studying degrees of real algebraic surfaces is $\R P^3$.
We may view a subset of $\R^3$ as a subset of $\R P^3$ which is
disjoint from the ``plane at infinity'', which is a copy of $\R P^2$.
We may define a \emph{flat} plane in $\R P^3$ to be the image of the
plane at infinity under a projective transformation of $\R P^3$, and a
\emph{topological} plane to be the image of the plane at infinity under
a homeomorphism of $\R P^3$.  This brings us to the following
generalization of the results of this paper:

\begin{conjecture} If a non-trivial link in $\R P^3$ is disjoint from
some topological plane, then it has four collinear points.
\end{conjecture}

\begin{conjecture} If an algebraic surface in $\R P^3$ is disjoint from some
topological plane and bounds a collection of non-trivially linked
solid tori, then the surface has degree at least 8.
\end{conjecture}

The following questions have also eluded the author:

\begin{conjecture} If an algebraic surface in $\R^3$ is a smooth torus
which is knotted on the outside, then it has degree at least eight.
\end{conjecture}

\begin{conjecture} Every wild arc in $\R^3$ has infinitely
many quadrisecants.
\end{conjecture}

\begin{question} What is the lowest possible degree of a polynomial
surface in $\R^3$ which is the boundary of the tubular neighborhood
of a trefoil knot?
\end{question}

% \bibliography{gt,books}

\begin{thebibliography}{1}

\bibitem{AA:geodesics}
Ralph Alexander and Stephanie Alexander, \emph{Geodesics in {Riemannian}
  manifolds with boundary}, Indiana Univ. Math. J. \textbf{30} (1981),
  481--488.

\bibitem{BZ:knots}
Gerhard Burde and Heiner Zieschang, \emph{Knots}, de Gruyter Studies in
  Mathematics, vol.~5, Walter de Gruyter, Berlin-New York, 1985.

\bibitem{Hartshorne:gtm}
Robin Hartshorne, \emph{Algebraic geometry}, Graduate Texts in Mathematics,
  vol.~52, Springer-Verlag, New York-Heidelberg-Berlin, 1977.

\bibitem{Moise:gtm}
Edwin~E. Moise, \emph{Geometric topology in dimensions 2 and 3}, Graduate Texts
  in Mathematics, vol.~47, Springer-Verlag, New York-Heidelberg-Berlin, 1977.

\bibitem{MM:quadrisecants}
H.~R. Morton and D.~M.~Q. Mond, \emph{Closed curves with no quadrisecants},
  Topology \textbf{21} (1982), 235--243.

\bibitem{Pannwitz:knoten}
E.~Pannwitz, \emph{Eine elmentargeometrische {Eigenschaft} von
  {Verschlingungen} und {Knoten}}, Math. Ann. \textbf{108} (1933), 629--672.

\bibitem{Reidemeister:knotentheorie}
K.~Reidemeister, \emph{Knotentheorie}, Chelsea Publishing Company, New York,
  1948.

\bibitem{Rolfsen:knots}
Dale Rolfsen, \emph{Knots and links}, Mathematics Lecture Series, vol.~7,
  Publish or Perish, Inc., Wilmington, DE, 1976.

\bibitem{Wall:generic}
C.~T.~C. Wall, Geometric properties of generic differentiable manifolds, Lect.
  Notes in Math., vol. 597, 707--774, Lect. Notes in Math., Springer-Verlag,
  New York-Heidelberg-Berlin, 1977, pp.~707--774.

\end{thebibliography}

\providecommand{\bysame}{\leavevmode\hbox to3em{\hrulefill}\thinspace}

\end{document}